# On the volume of projections of the cross-polytope

*Grigory Ivanov*[1]

**Abstract.** We study properties of the volume of projections of the $n$-dimensional cross-polytope $\Diamond^n = \{x \in \mathbb{R}^n \mid |x_1| + \cdots + |x_n| \leqslant 1\}$. We prove that the projection of $\Diamond^n$ onto a $k$-dimensional coordinate subspace has the maximum possible volume for $k = 2$ and for $k = 3$. We obtain the exact lower bound on the volume of such a projection onto a two-dimensional plane. Also, we show that there exist local maxima which are not global ones for the volume of a projection of $\Diamond^n$ onto a $k$-dimensional subspace for any $n > k \geqslant 2$.

## 1. Introduction

The standard cross-polytope $\Diamond^n$ in $\mathbb{R}^n$ is the convex hull of the vectors of the standard basis $\{e_i\}_1^n$ of $\mathbb{R}^n$ and their opposites, $\{-e_i\}_1^n$. The cross-polytope is the dual (polar) body for the standard cube $\square^n = [-1,1]^n$. In this paper we study lower and upper bounds on the volume of a projection of $\Diamond^n$ onto a $k$-dimensional subspace, that is, we study extremizers of

$$(1.1) \qquad F(H_k) = \operatorname{vol}_k\left(\Diamond^n | H_k\right),$$

where $H_k$ is a $k$-dimensional subspace of $\mathbb{R}^n$ and $\Diamond^n | H_k$ denotes the projection of $\Diamond^n$ onto $H_k$. This problem can be considered as the dual problem for finding extremizers of

$$F_2(H_k) = \operatorname{vol}_k\left(\square^n \cap H_k\right),$$

whose global extrema are reasonably well-studied

$$(1.2) \qquad \operatorname{vol}_k \square^k \leqslant \operatorname{vol}_k\left(\square^n \cap H_k\right) \leqslant \min\left\{\left(\frac{n}{k}\right)^{\frac{k}{2}}, 2^{\frac{n-k}{2}}\right\} \operatorname{vol}_k \square^k,$$

where the left-hand side inequality is due to J. Vaaler [Vaa79] and it is optimal for all $n \geqslant k \geqslant 1$; and the right-hand side is due to K. Ball [Bal89] and it is optimal whenever $k|n$ or $2k \geqslant n$ (that is, the optimal constants are not known only when $2k < n$ and $k \nmid n$). To get a rather complete survey on the volumes of intersections of the cube with linear subspaces, we refer the reader to the first chapter of Zong's book [Zon06]. A very interesting generalization of the Vaaler result was made by M. Meyer and A. Pajor [MP88] and a new proof of Vaaler's inequality in terms of waists was given in [AHK19] recently.

Surprisingly, the optimal constants in $C_1 \operatorname{vol}_k \Diamond^k \leqslant F(H_k) \leqslant C_2 \operatorname{vol}_k \Diamond^k$ are mostly not known. The following conjecture can be considered as an analogue (or the dual conjecture) for Vaaler's inequality.

**Conjecture 1.1.** *The volume of the projection of the $n$-dimensional cross-polytope $\Diamond^n$ onto a $k$-dimensional subspace is at most the volume of the $k$-dimensional cross-polytope, i.e. $\operatorname{vol}_k(\Diamond^n | H_k) \leqslant \frac{2^k}{k!}$. The bound is attained only on coordinate subspaces.*

The only known results here are that the conjecture is true in the hyperplane case $k = n - 1$ (see [Bal95]), which is quite simple (since the projections of the surface of the cross-polytope onto a hyperplane covers the projection of the cross-polytope twice almost everywhere) and in the two-dimensional case [Fil88]. Also, P. Filliman proved that this conjecture is true for $k = 3$ and $n \leqslant 6$. The deep generalization of the hyperplane case for the volume of a projection of the $\ell_p^n$ balls was given by A. Naor and F. Barthe [BN02]. We show that this Conjecture is true in lower-dimensional cases.

[1] Institute of Discrete Mathematics and Geometry, TU Wien, Vienna.
Department of Higher Mathematics, Moscow Institute of Physics and Technology, Institutskii pereulok 9, Dolgoprudny, Moscow region, 141700, Russia.
grimivanov@gmail.com
Research was supported by the Swiss National Science Foundation grant 200021-179133 and by the Russian Foundation for Basic Research, project 18-01-00036A.

**Theorem 1.2.** *Conjecture 1.1 is true for $k = 2, 3$ and any $n \geq k$.*

To prove this theorem, we give in Theorem 5.3 a geometric first-order necessary condition for $\lozenge^n | H_k$ (for any $k$ and $n \geq k$) to be a maximizer of (1.1). In Lemma 4.8, we show that for a fixed $k$ and an arbitrary $n \geq k$ the maximum of (1.1) is already attained on some $n \leq k^3$. This implies that there are finitely many combinatorial structures that the extremizer $\lozenge^n | H_k$ might have, and the maximum can be found by computation for a small $k$.

The only known optimal subspaces $H_k$ in Ball's inequality, which is the second inequality in (1.2), are such that the section of the cube by $H_k$ is an affine cube. In [Iva17], the author proved that all optimal subspaces have this property whenever $k \mid n$. We think that this is the case for an arbitrary $k \geq 2$ and $n \geq k$. Let $C_\square(n,k) 2^k$ be the maximum volume of a section of $\square^n$ by a $k$-dimensional linear subspace $H_k$ such that $\square^n \cap H_k$ is an affine cube.

**Conjecture 1.3.** *The maximum volume of a section of $\square^n$ by a $k$-dimensional linear subspace $H_k$ is attained on subspaces such that $\square^n \cap H_k$ is an affine cube, i.e. $\operatorname{vol}_k(\square^n \cap H_k) \leq C_\square(n,k) 2^k$.*

And we can formulate the dual statements for projections of $\lozenge^n$. Let $C_\lozenge(n,k) \frac{2^k}{k!}$ be the minimum volume of a projection of the $\lozenge^n$ onto a $k$-dimensional subspace $H_k$ such that $\lozenge^n | H_k$ is an affine cross-polytope.

**Conjecture 1.4.** *The minimum volume of a projection of $\lozenge^n$ onto a $k$-dimensional subspace $H_k$ is attained on subspaces such that $\lozenge^n | H_k$ is an affine cross-polytope, i.e. $\operatorname{vol}_k(\lozenge^n | H_k) \geqslant C_\lozenge(n,k) \frac{2^k}{k!}$.*

It is not hard to prove that

$$(1.3) \qquad C_\square^2(n,k) = \frac{1}{C_\lozenge^2(n,k)} = \left(\left\lceil \frac{n}{k} \right\rceil\right)^{n-k\lfloor \frac{n}{k} \rfloor} \left(\left\lfloor \frac{n}{k} \right\rfloor\right)^{k-\left(n-k\lfloor \frac{n}{k} \rfloor\right)},$$

and that these constants are attained on the same subspaces. For completeness, we give a proof of (1.3) and give a complete description of a $k$-dimensional subspace $H_k$ on which constants $C_\square(n,k)$ and $C_\lozenge(n,k)$ are attained in Lemma 4.7.

K. Ball proved that the constant in Conjecture 1.3 is optimal whenever $k|n$ or $2k \geqslant n$, F. Barthe [Bar98] proved the constant in Conjecture 1.4 is optimal whenever $k|n$ and, due to F. Barthe and A. Naor [BN02], it is optimal in the hyperplane case $k = n - 1$. In [Iva17, Theorem 1.5], the author proved that the last two conjectures are true in case $k|n$. We can give only a partial answer to Conjecture 1.4 in the most simple case.

**Theorem 1.5.** *Conjecture 1.4 is true for $k = 2$ and any $n \geq k$.*

It is well-known that any centrally symmetric polytope is an affine image of the high-dimensional cross-polytope or an affine image of the intersection of the high-dimensional cube by a linear subspace. This means that we can reformulate our conjectures in terms of uniform bounds on the volume of a centrally symmetric polytope in $\mathbb{R}^k$. We formulate the proper inequalities in Section 3 after some preliminaries and explanations.

The paper is organized as follows: After presenting our notation and basic definitions, we reformulate the problems in terms of the polytopes generated by projections of vectors of an orthonormal basis in Section 2. Then in Section 3, following the ideas of [Iva18], we explain how to obtain a first order necessary condition for extremizers. In Section 4, we simplify the structure of local extremizers of (1.1). Then, after some geometric constructions, we give a geometrical necessary condition for $H_k$ to be a local maximizer of (1.1). Using these results, we prove Theorems 1.5 and 1.2 in Section 6. Then, in Section 7, we show that there exist local and non-global maximizers of (1.1).

## 2. Defintions and Preliminaries

For a positive integer $n$, we refer to the set $\{1, 2, \ldots, n\}$ as $[n]$. We use $\Diamond^n$ to denote an $n$-dimensional cross-polytope $\{x | \sum_{i \in [n]} |x[i]| \leqslant 1\}$ in $\mathbb{R}^n$. Here and throughout the paper $x[i]$ stands for the $i$-th coordinate of a vector $x$. As usual, $\{e_i\}_1^n$ is the standard orthonormal basis of $\mathbb{R}^n$. We use $\langle p, x \rangle$ to denote the *value of a linear functional $p$ at a vector $x$*. For vectors $u, v \in \mathbb{R}^d$, their *tensor product* (or, diadic product) is a linear operator on $\mathbb{R}^d$ defined as $(u \otimes v)x = \langle u, x \rangle v$ for every $x \in \mathbb{R}^d$.

Throughout the paper $H_k$ will be a $k$-dimensional subspace of $\mathbb{R}^n$. For a convex body $K \subset \mathbb{R}^n$ and a $k$-dimensional subspace $H_k$ of $\mathbb{R}^n$ we denote by $K \cap H_k$ and $K|H_k$ the section of $K$ by $H_k$ and the orthogonal projection of $K$ onto $H_k$, respectively. For a $k$-dimensional subspace $H_k$ of $\mathbb{R}^n$ and a convex body $K \subset H_k$ we denote by $\mathrm{vol}_k K$ the $k$-dimensional volume of $K$. We consider only $n \geqslant k \geqslant 2$.

First of all, it is convenient to identify a projection of the cross-polytope with a convex polytope in $\mathbb{R}^k$. Let $v_i = Pe_i$, where $P$ is the projection onto $H_k$. Clearly,
$$\Diamond^n | H_k = \mathrm{co}\{\pm v_1, \ldots, \pm v_n\}.$$
That means that a projection of $\Diamond^n$ is determined by the set of vectors $(v_i)_1^n$, which are the projections of the orthogonal basis. Such sets of vectors have several equivalent description and names.

**Definition 2.1.** We will say that an ordered $n$-tuple of vectors $\{v_1, \ldots, v_n\} \subset H$ forms *a tight frame* in a vector space $H$ if

(2.1) $$\left( \sum_1^n v_i \otimes v_i \right) \bigg|_H = I_H,$$

where $I_H$ is the identity operator in $H$ and $A|_H$ is the restriction of an operator $A$ onto $H$. We use $\Omega(n, k)$ to denote the set of all tight frames with $n$ vectors in $\mathbb{R}^k$.

For the sake of convenience, we write $I_k$ for $I_{\mathbb{R}^k}$.

**Definition 2.2.** An ordered $n$-tuple of vectors in a linear space $H$ which span $H$ is called a *frame*.

In the following trivial Lemma we understand $\mathbb{R}^k \subset \mathbb{R}^n$ as the subspace of vectors, whose last $n - k$ coordinates are zero. For convenience, we will consider $\{v_i\}_1^n \subset \mathbb{R}^k \subset \mathbb{R}^n$ to be $k$-dimensional vectors.

**Lemma 2.3.** *The following assertions are equivalent:*
  *(1) the vectors $\{v_1, \ldots, v_n\} \subset \mathbb{R}^k$ form a tight frame in $\mathbb{R}^k$;*
  *(2) there exists an orthonormal basis $\{f_1, \ldots, f_n\}$ of $\mathbb{R}^n$ such that $v_i$ is the orthogonal projection of $f_i$ onto $\mathbb{R}^k$, for any $i \in [n]$;*
  *(3) $\mathrm{Lin}\{v_1, \ldots, v_n\} = \mathbb{R}^k$ and the Gram matrix $\Gamma$ of vectors $\{v_1, \ldots, v_n\} \subset \mathbb{R}^k$ is the matrix of a projection operator from $\mathbb{R}^n$ onto the linear hull of the rows of the matrix $M = (v_1, \ldots, v_n)$.*
  *(4) the $k \times n$ matrix $M = (v_1, \ldots, v_n)$ is a sub-matrix of an orthogonal matrix of order $n$.*

It follows that the tight frames in $\mathbb{R}^k$ are exactly the projections of orthonormal bases onto $\mathbb{R}^k$. This observation allows to reformulate the problems in terms of tight frames and associated polytopes in $\mathbb{R}^k$. Indeed, since $H_k$ and $\mathbb{R}^k$ are isometric, $\Diamond^n | H_k$ can be considered as the absolute convex hull of a tight frame formed by the projections of the standard basis onto $H_k$. Vice versa, assertion (3) gives a way to reconstruct $H_k$ from a given tight frame $\{v_1, \ldots, v_n\}$ in $\mathbb{R}^k$.

**Definition 2.4.** We will say that the set $S = \{v_1, \ldots, v_n\}$ of $n$ vectors of $\mathbb{R}^k$ generates

(1) a polytope $\mathrm{co}\{\pm v_1,\ldots,\pm v_n\}$ We use $\lozenge^n|S$ to denote this polytope. We will say that $\lozenge^n|S$ is the *associated with $S$ polytope*.
(2) a matrix $\sum_{i\in[n]} v_i \otimes v_i$. We use $A_S$ to denote this matrix.

Summarizing, the global extrema of (1.1) coincide with those of

(2.2) $$F(S) = \mathrm{vol}_k(\lozenge^n|S), \qquad \text{where} \qquad S \in \Omega(n,k).$$

In order to compare local extrema of both problems and write the first order necessary condition, we need to endow the Grassmannian of $k$-dimensional subspaces of $\mathbb{R}^n$ and the set $\Omega(n,k)$ with metrics. The standard metric on Grassmannian is given by

$$\mathrm{Dist}(H,H') = \left\|P^H - P^{H'}\right\|,$$

where $\|\cdot\|$ denotes the operator norm, and $P^H$ and $P^{H'}$ are the orthogonal projections onto subspaces $H$ and $H'$, respectively. We endow $\Omega(n,k)$ with a metric

$$\mathrm{dist}(\{v_1,\ldots,v_n\},\{w_1,\ldots,w_n\}) = \sqrt{\sum_1^n |v_i - w_i|^2}.$$

As was shown by the author in [Iva18, Section 4], $H_k$ is a local extremizer of (1.1) iff a tight frame formed by the projections of the standard basis of $\mathbb{R}^n$ is a local extremizer of (2.2) (for any isometry between $H_k$ and $\mathbb{R}^k$). This and an observation that any frame can be mapped to a tight frame allow us to consider small perturbations of a tight frame $S \in \Omega(n,k)$ and associated polytope $\lozenge^n|S$. The advantage is that one can understand what happens geometrically with the polytope $\lozenge^n|S$ and its volume after perturbations of its vertices.

For the sake of completeness, we explain the approach proposed in [Iva18] of obtaining the first-order necessary condition for (2.2) in the next Section.

## 3. Operations on frames

Our main idea is to transform a given tight frame $S$ to a new one $S'$ and compare the volumes of the projections of the cross-polytope generated by them. Since it is not very convenient to transform a given tight frame to another tight frame, we add an intermediate step – we transform a tight frame $S$ to a frame $\tilde{S}$, and then we transform $\tilde{S}$ to a new tight frame $S'$ using a linear transformation. The main observation here is that we can always transform any frame $\tilde{S} = \{v_1,\ldots,v_n\}$ to a tight frame $S'$ using a proper linear transformation $L$: $S' = L\tilde{S} = \{Lv_1,\ldots,Lv_n\}$, or, equivalently, any non-degenerate centrally symmetric polytope in $\mathbb{R}^k$ is an affine image of a projection of a high dimension cross-polytope.

For a frame $S = \{v_1,\ldots,v_n\}$ in $\mathbb{R}^k$, by definition put

$$B_S = A_S^{-\frac{1}{2}} = \left(\sum_{i\in[n]} v_i \otimes v_i\right)^{-\frac{1}{2}}.$$

The operator $B_S$ is well-defined as the condition $\mathrm{Lin}\, S = \mathbb{R}^k$ implies that $A_S$ is a positive definite operator. Clearly, $B_S$ maps any frame $S$ to a tight frame:

$$\sum_{i=1}^n B_S v_i \otimes B_S v_i = B_S \left(\sum_{i=1}^n v_i \otimes v_i\right) B_S^T = B_S A_S B_S = I_k.$$

This means that Conjectures 1.1, 1.3 and 1.4 can be rewritten in the following way.

**Conjecture 3.1.** *Let $S = \{v_1,\ldots,v_n\}$ be a frame in $\mathbb{R}^k$ and by definition put $P = \lozenge^n|S$ and $Q = \cap_1^n \{x \in \mathbb{R}^k \mid |\langle x, v_i\rangle| \leqslant 1\}$. Then*

$$C_\lozenge(n,k) \mathrm{vol}_k \lozenge^k \leqslant \mathrm{vol}_k P \cdot \det B_S \leqslant \mathrm{vol}_k \lozenge^k; \quad \mathrm{vol}_k \square^k \leqslant \frac{\mathrm{vol}_k Q}{\det B_S} \leqslant C_\square(n,k) \mathrm{vol}_k \square^k.$$

We obtain the following necessary and sufficient condition for a tight frame to be a maximizer of (2.2).

**Lemma 3.2.** *The maximum (respectively minimum) in* (2.2) *is attained on a tight frame* $S = \{v_1, \ldots, v_n\}$ *iff the following inequality holds for an arbitrary frame* $\tilde{S}$

$$\text{(3.1)} \qquad \frac{\text{vol}_k(\lozenge^n|\tilde{S})}{\text{vol}_k(\lozenge^n|S)} \leqslant \sqrt{\det A_{\tilde{S}}} \quad (\text{respectively} \geqslant).$$

**Proof.**
As mentioned above, $B_{\tilde{S}}\tilde{S}$ is a tight frame, and clearly, $\text{vol}_k(\lozenge^n|B_{\tilde{S}}\tilde{S}) = \det B_{\tilde{S}} \, \text{vol}_k(\lozenge^n|\tilde{S})$. The tight frame $S$ is a maximizer iff $\text{vol}_k(\lozenge^n|B_{\tilde{S}}\tilde{S}) \leqslant \text{vol}_k(\lozenge^n|S)$ for an arbitrary frame $\tilde{S}$. Using these observations and the definition of $B_S$, we have

$$\frac{\text{vol}_k(\lozenge^n|\tilde{S})}{\text{vol}_k(\lozenge^n|S)} = \frac{1}{\det B_{\tilde{S}}} \frac{\text{vol}_k(\lozenge^n|B_{\tilde{S}}\tilde{S})}{\text{vol}_k(\lozenge^n|S)} \leqslant \frac{1}{\det B_{\tilde{S}}} = \sqrt{\det A_{\tilde{S}}}.$$

□

Let us illustrate how we will use it. Let $S$ be an extremizer of (2.2), and $T$ be a map from a subset of $\Omega(n,k)$ to the set of frames. In order to obtain properties of extremizers, we consider a composition of two operations:

$$\text{(3.2)} \qquad S \xrightarrow{\mathbf{T}} \tilde{S} \xrightarrow{B_{\tilde{S}}} S',$$

where $B_{\tilde{S}}$ is as defined above. For example, see Figure 1, where $T$ is the operation of replacing a vector $v$ of $S$ by the origin.

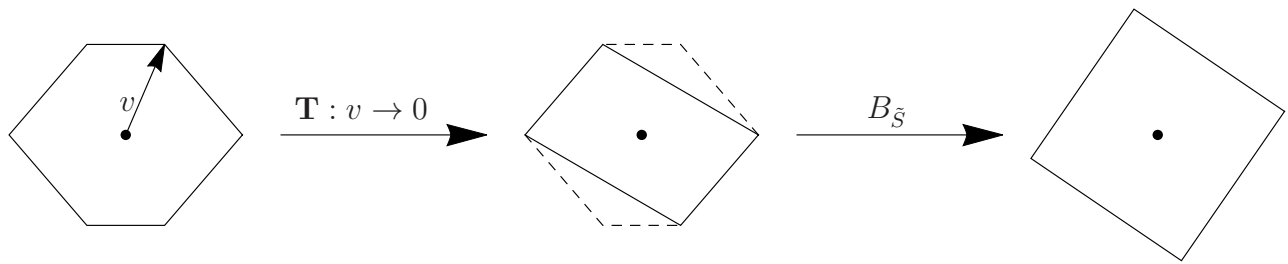

FIGURE 1. Here we map one vector to zero.

Choosing a simple operation $\mathbf{T}$, we may calculate the left-hand side of (3.1) in some geometric terms. We consider several simple operations: Scaling one or several vectors, mapping one vector to the origin, mapping one vector to another. On the other hand, the determinant in the right-hand side of (3.1) can be calculated for the operations listed above. In particular, the following first-order approximation of the determinant was obtained by the author in [Iva18, Theorem 1.2]. We provide a sketch of its proof in the Appendix 8.

**Lemma 3.3.** *For an arbitrary tight frame $S$ the following identity holds*

$$\sqrt{\det A_{S'}} = 1 + \sum_{i=1}^{n} t_i \langle x_i, v_i \rangle + o\left(\sqrt{t_1^2 + \cdots + t_n^2}\right),$$

*where $S'$ is obtained from $S$ by substitution $v_i \to v_i + t_i x_i$, $i \in [n]$.*

Clearly, if $\tilde{S}$ in the assertion of Lemma 3.2 is close to $S$, then the tight frame $B_{\tilde{S}}\tilde{S}$ is close to $S$ as well. Therefore, inequality (3.1) gives necessary condition for local maximizers of (2.2).

Also, since the composition of projections from $\mathbb{R}^n$ onto $\mathbb{R}^k$ and from $\mathbb{R}^k$ onto a subspace $H \subset \mathbb{R}^k$ is a projection, we have the following.

**Lemma 3.4.** *Let $H$ be a subspace of $\mathbb{R}^k$ and $P$ be the orthogonal projection onto $H$, let $S = \{v_1, \ldots, v_n\}$ be a tight frame in $\mathbb{R}^k$. Then the vectors $\{Pv_1, \ldots, Pv_n\}$ form a tight frame in $H$.*

## 4. Reduction of the problems

We need some simple properties of the determinants of positive definite operators.

**Lemma 4.1.** *Let $M$ be an operator $\mathbb{R}^d \to \mathbb{R}^d$ given by $I_d \pm v \otimes v$, where $I_d$ is the identity operator on $\mathbb{R}^d$ and $v \in \mathbb{R}^d$. Then*

(4.1) $$\det M = 1 \pm |v|^2.$$

**Proof.**
The lemma is trivial. Operator $M$ is diagonal in an orthonormal basis of $\mathbb{R}^d$ in which the first basis vector is collinear to $v$. The first entry is $1 \pm |v|^2$. All others are equal to one.
□

Since $\operatorname{tr} A_S = \sum_1^n \operatorname{tr} v_i \otimes v_i = \sum_1^n |v_i|^2$ for any $n$-tuple of vectors of $\mathbb{R}^d$ and by definition of a tight frame, we get the following trivial property.

**Lemma 4.2.** *Let $S = \{v_1, \ldots, v_n\}$ be a tight frame in $\mathbb{R}^k$. Then*

(4.2) $$\operatorname{tr} A_S = \sum_1^n |v_i|^2 = k.$$

Let $S \setminus i$ denote the set $(v_j)_{j \in [n]}^{j \neq i}$ for a frame $S = \{v_1, \ldots, v_n\}$,

**Lemma 4.3.** *Let $S'$ be a set of vectors obtained from a tight frame $S = \{v_1, \ldots, v_n\}$ by substitution $v_i \to v_i + tx$, where $t \in \mathbb{R}$, $x \in \mathbb{R}^k$ and a fixed index $i$. If $|v_i| < 1$, then*

(4.3) $$\det B_{S'} = \sqrt{\frac{1 + |B_{S \setminus i} v_i|^2}{1 + |B_{S \setminus i} v|^2}},$$

*where $v = v_i + tx$.*

**Proof.**
Since $|v_i| < 1$, we have that $A_{S \setminus i} = I - v_i \otimes v_i$ is a positive definite matrix. Hence, $B_{S \setminus i} > 0$.

Let $\tilde{S} = \{v'_1, \ldots, v'_n\}$ be the frame obtained from $S'$ by applying the linear transformation $B_{S \setminus i}$. That is, $v'_j = B_{S \setminus i} v_j, j \in [n] \setminus i$, and $v'_i = B_{S \setminus i} v$. By the definition of $B_S$, we have that $A_{\tilde{S} \setminus i} = I_k$ and $A_{\tilde{S}} = I_k + v'_i \otimes v'_i$. Hence, by (4.1), we get $\det A_{\tilde{S}} = 1 + |v'_i|^2 = 1 + |B_{S \setminus i} v|^2$.

By the same argument,
$$\det A_{B_{S \setminus i} S} = 1 + |B_{S \setminus i} v_i|^2.$$

As $\det A_S = 1$, we obtain

$$\det B_{S'} = \frac{1}{\sqrt{\det A_{S'}}} = \sqrt{\frac{\det A_S}{\det A_{S'}}} = \sqrt{\frac{\det A_{B_{S \setminus i} S}}{\det A_{B_{S \setminus i} S'}}} = \sqrt{\frac{1 + |B_{S \setminus i} v_i|^2}{1 + |B_{S \setminus i} v|^2}}.$$

□

Using results from the previous section, we can simplify a configuration of extremizers of (2.2).

The straightforward consequence of the following lemma is that if $S = \{v_1, \ldots, v_n\}$ is a (local) maximizer of (2.2) then the non-zero elements of $\{\pm v_1, \ldots, \pm v_n\}$ are pairwise distinct vertices of $\lozenge^n | S$.

**Lemma 4.4.** *Let $S = \{v_1, \ldots, v_n\}$ be a tight frame such that $v_i \in \Diamond^{n-1}|(S \setminus i)$. Let $\tilde{S}$ be the frame obtained from $S$ by substitution $v_i \to 0$. By definition put $S' = B_{\tilde{S}}\tilde{S}$. Then*
$$\mathrm{vol}_k(\Diamond^n|S) \leqslant \mathrm{vol}_k(\Diamond^n|S'),$$
*and equality holds iff $v_i = 0$.*

**Proof.**
Indeed, we have that $\Diamond^n|S = \Diamond^n|\tilde{S}$ and, obviously, $B_{\tilde{S}} \geqslant B_S = I$ (moreover, the equality holds iff $v_i = 0$). Therefore,
$$\mathrm{vol}_k(\Diamond^n|S') = \det B_{\tilde{S}} \, \mathrm{vol}_k(\Diamond^n|\tilde{S}) \geqslant \mathrm{vol}_k(\Diamond^n|\tilde{S}) = \mathrm{vol}_k(\Diamond^n|S).$$
□

Moreover, we can change $v_i$ in a continuous way $v_i \to (1-\lambda)v_i$, where $\lambda \in [0,1]$, while increasing the volume of $\Diamond^n|S'$.

Using the same arguments and (4.3), we get

**Corollary 4.5.** *Let $\tilde{S}$ be the frame obtained from $S$ by substitution $v_i \to v$. If $|v_i| < 1$, $|B_{S\setminus i}v_i| \leqslant |B_{S\setminus i}v|$ and $\Diamond^n|\tilde{S} \subset \Diamond^n|S$, then for the tight frame $S' = B_{\tilde{S}}\tilde{S}$, we have*
$$\mathrm{vol}_k(\Diamond^n|S') \leqslant \mathrm{vol}_k(\Diamond^n|S).$$
*Equality holds iff $|B_{S\setminus i}v_i| = |B_{S\setminus i}v|$ and $\Diamond^n|\tilde{S} = \Diamond^n|S$.*

Now we are ready to prove some simple properties of extremizers of (2.2). Let us start with the minimizers.

**Lemma 4.6.** *If $\min\limits_{S \in \Omega(n,k)} \mathrm{vol}_k(\Diamond^n|S)$ is attained on $S = \{v_1, \ldots, v_n\}$, then*

*(1) all vectors $|v_i|$ are vertices of $\Diamond^n|S$;*
*(2) $|v_j| \leqslant \sqrt{2}|v_i|$, for any $i, j \in [n]$;*
*(3) at most $\lfloor \frac{n}{k} \rfloor + 1$ vectors of $S$ may be collinear.*

**Proof.**
1) Suppose vector $v_i$ is not a vertex of $\Diamond^n|S$. Then $v_i$ is a convex combination of some other vectors of $S$ and $-S$. The same is true for $B_{S\setminus i}v_i$ and $B_{S\setminus i}S$. Therefore, by the triangle inequality, we have that $|v_i| < 1$ and there exists $j \in [n]$ such that $|B_{S\setminus i}v_i| < |B_{S\setminus i}v_j|$. Let $\tilde{S}$ be the frame obtained from $S$ by substitution $v_i \to v_j$. We have that $|v_i| < 1$, $|B_{S\setminus i}v_i| < |B_{S\setminus i}v|$ and $\Diamond^n|\tilde{S} \subset \Diamond^n|S$, Applying Corollary 4.5 for the tight frame $S$, we get a tight frame $S'$ such that $\mathrm{vol}_k(\Diamond^n|S') < \mathrm{vol}_k(\Diamond^n|S)$. This contradicts to the initial choice of $S$.

2) If $|v_i| \geqslant \frac{1}{\sqrt{2}}$, there is nothing to prove. Assume $|v_i| < \frac{1}{\sqrt{2}}$ and $|v_j| > \sqrt{2}|v_i|$ for a fixed $j \in [n]$. Let $\tilde{S}$ be the frame obtained from $S$ by substitution $v_i \to v_j$.

Clearly, $0 < A_{S\setminus i} \leqslant A_S = I_k$ and, consequently, $B_{S\setminus i} \geqslant I_k$. Therefore, $|B_{S\setminus i}v_j| \geqslant |v_j|$. Since $A_{S\setminus i} = I_k - v_i \otimes v_i$ scales the space in the direction $v_i$ with factor $1 - |v_i|^2$ and it is the identity on the hyperplane $v_i^\perp$, we have that $B_{S\setminus i} = (A_{S\setminus i})^{-1/2}$ scales the space in the same direction with factor $\frac{1}{\sqrt{1-|v_i|^2}}$. That is, $B_{S\setminus i}v_i = \frac{v_i}{\sqrt{1-|v_i|^2}}$. One can see that
$$\left|\frac{v_i}{\sqrt{1-|v_i|^2}}\right| \leqslant \sqrt{2}|v_i| < |v_j| \leqslant |B_{S\setminus i}v_j|.$$

Using the same arguments as above, we get that $\mathrm{vol}_k(\Diamond^n|S') < \mathrm{vol}_k(\Diamond^n|S)$, where $S' = B_{\tilde{S}}\tilde{S}$ This completes the proof of the second assertion of the lemma.

3) Assume the contrary. Let the vectors $v_1, \ldots, v_d$ be collinear, where $d \geqslant \lfloor \frac{n}{k} \rfloor + 2$. By the assertion (1), we can assume that $v_1 = \cdots = v_d$. Consider the affine transformation $B_{S\setminus 1}$. We have that $B_{S\setminus 1}v_1 = \cdots = B_{S\setminus 1}v_d$ and that the vectors $B_{S\setminus 1}v_2, \ldots, B_{S\setminus 1}v_n$ form a tight frame in $\mathbb{R}^k$. Applying Lemma 2.3, we get

$$|B_{S\setminus 1}v_2|^2 + \cdots + |B_{S\setminus 1}v_d|^2 \leqslant 1.$$

Therefore, $|B_{S\setminus 1}v_1| = |B_{S\setminus 1}v_2| \leqslant \sqrt{\frac{1}{d-1}}$. Again, since the vectors $B_{S\setminus 1}v_2, \ldots, B_{S\setminus 1}v_n$ form a tight frame and by Lemma 2.3, there exists $j \in [n]$ such that

$$|B_{S\setminus 1}v_j| \geqslant \sqrt{\frac{k}{n-1}} > \sqrt{\frac{k}{n}} > \sqrt{\frac{1}{d-1}} \geqslant |B_{S\setminus 1}v_1|.$$

Let $\tilde{S}$ be the frame obtained from $S$ by substitution $v_1 \to v_j$ and $S' = B_{\tilde{S}}\tilde{S}$. Using the same arguments as above, we get that $\operatorname{vol}_k(\Diamond^n|S') < \operatorname{vol}_k(\Diamond^n|S)$. This contradicts the initial choice of $S$.
$\square$

And again, in all assertions of this lemma we can transform $S$ to $S'$ in a continuous way while decreasing the volume in a monotonic way. And now, we are ready to describe all $H_k$ such that $C_\Diamond(n,k)$ is attained.

**Lemma 4.7.** *The constants $C_\Diamond(n,k)$ and $C_\Box(n,k)$ are given by (1.3) and they are attained on the same subspaces. All such subpspaces are given by the following rule:*

*(1) we partition $[n]$ into $k$ sets such that the cardinalities of any two sets differ by at most one;*

*(2) let $\{i_1, \ldots, i_\ell\}$ be one of the sets of the partition. Then, choosing arbitrary signs, we write the system of linear equations*

$$\pm x[i_1] = \cdots = \pm x[i_\ell];$$

*(3) our subspace is the solution of the system of all equations written for each set of the partition at the step (2).*

**Proof.**
By a duality argument, $\Diamond^n|H_k$ is an affine cross-polytope iff $\Box^n \cap H_k$ is an affine cube and they are polars in $H_k$. Therefore, $\operatorname{vol}_k(\Diamond^n|H_k) \cdot \operatorname{vol}_k(\Box^n \cap H_k) = \operatorname{vol}_k\Diamond^k \cdot \operatorname{vol}_k\Box^k$ in case $\Diamond^n|H_k$ is an affine cross-polytope. Hence, the first equation (1.3) is proven.

Let $H_k$ be such that $C_\Diamond(n,k)$ is attained and let $\{\pm a_1, \ldots, \pm a_k\}$ be the vertices of $\Diamond^n|H_k$. By the first assertion of Lemma 4.6, we have that the projections $v_i$ of the vectors $e_i$ are the vertices of $\Diamond^n|H_k$, i.e. $v_i$ coincides with $\pm a_j$ for a proper sign and $j \in [k]$. Or, equivalently, we partition $[n]$ into $k$ sets and $H_k$ is the solution of a proper system of linear equations constructed as in (2) and (3), except we have not proved that (1) holds yet. Let us prove this assertion. Let $d_i$ vectors of the standard basis of $\mathbb{R}^n$ project onto a pair $\pm a_i$. Therefore, a $k$-tuple of vectors $\{\sqrt{d_i}a_i\}_{i\in[k]}$ forms a tight frame. Identifying $H_k$ with $\mathbb{R}^k$ and by the assertion (3) of Lemma 2.3, we conclude that $a_i$ and $a_j$ are orthogonal whenever $i \neq j$. Therefore, $|a_i|^2 = \frac{1}{d_i}$ and

$$(4.4) \qquad \operatorname{vol}_k(\Diamond^n|H_k) = \frac{2^k}{k!} \frac{1}{\sqrt{d_1 \cdot \ldots \cdot d_k}}.$$

Suppose $d_i \geqslant d_j + 2$ for some $i, j \in [k]$. Then $d_i \cdot d_j \leqslant (d_i - 1)(d_j + 1)$. By this and by (4.4), we showed that (1) holds.

It is easy to see that there are exactly $n - k\lfloor\frac{n}{k}\rfloor$ of $d_i$'s equal $\lceil\frac{n}{k}\rceil$ and all others $k - (n - k\lfloor\frac{n}{k}\rfloor)$ are equal to $\lfloor\frac{n}{k}\rfloor$. That is, $C_\Diamond(n,k)$ is given by (1.3). This completes the proof.
$\square$

The simple Lemma 4.4 allows us to bound efficiently $n$ for a maximizer of (2.2).

**Lemma 4.8.** *Denote $V(n,k) = \max_{S \in \Omega(n,k)} \operatorname{vol}_k(\Diamond^n|S)$. Then $V(n,k) \leqslant V(k^3, k)$.*

**Proof.**
Let us show that $V(n,k) \geqslant \frac{2}{k} V(n-1, k-1)$. Indeed, consider a tight frame $S_1$ of $n-1$ vectors in $\mathbb{R}^{k-1}$ such that $\mathrm{vol}_{k-1}\left(\Diamond^{n-1}|S_1\right) = V(n-1, k-1)$. Adding $e_k$ to $S_1$, we obtain a new tight frame, lets call it $S$. Obviously, $\mathrm{vol}_k\left(\Diamond^n|S\right) = \frac{2}{k} \mathrm{vol}_{k-1}\left(\Diamond^{n-1}|S_1\right)$. The needed inequality is proven.

Let $S = \{v_1, \ldots, v_n\}$ be a tight frame in $\mathbb{R}^k$ such that $\mathrm{vol}_k\left(\Diamond^n|S\right) = V(n,k)$, and let $n > k^3$. By Lemma 4.4, we conclude that all vectors of $S$ are non-zero vectors (otherwise we can omit them and decrease $n$). Since $\sum_{1}^{n} |v_i|^2 = k$, there exists $i$ such that $|v_i| \leqslant \sqrt{\frac{k}{n}}$. By Lemma 3.4, the projections of $\{v_1, \ldots, v_n\}$ onto the hyperplane $H_i$ perpendicular to $v_i$ form a tight frame in $H_i$. As $v_i$ projects at zero, we may consider all others projections as a tight frame in $H_i$, which we denote $S_1$. On the other hand, we have

$$V(n,k) = \mathrm{vol}_k\left(\Diamond^n|S\right) \leqslant 2|v_i| \mathrm{vol}_{k-1}\left(\Diamond^{n-1}|S_1\right) < \frac{2}{k} V(n-1, k-1).$$

We come to a contradiction.
□

## 5. Geometric constructions and properties of maximizers

As mentioned above, to use Lemma 3.2, we need to understand the geometry behind the left-hand side of (3.1). For this purpose we introduce the following definitions.

Let $\{\pm v_1, \ldots, \pm v_n\}$ be a set of pairwise distinct vectors and be the vertex set of a centrally symmetric polytope $P$.

We use $\mathcal{F}(v)$ to denote the set of facets of $P$ incident to its vertex $v$. Then we define the *star* of the vertex $v$ of $P$ as follows

$$N_P(v) = \bigcup_{F \in \mathcal{F}(v), F \in \mathcal{F}(-v)} \mathrm{co}\{F, 0\}.$$

By definition put $Q_P(v) = \overline{P \setminus N_P(v)}$ and $R_P(v) = \mathrm{co}\{\pm v_i | v_i \neq \pm v\}$, we call these sets the *belt* and the *rest* of the vertex $v$ in $P$, respectively.

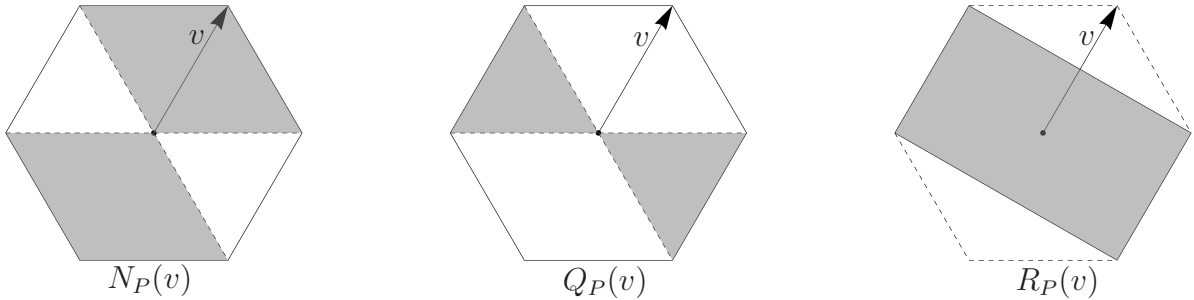

FIGURE 2. The star, the belt and the rest of vertex $v$ for the regular hexagon $P$.

By symmetry, $N_P(v) = N_P(-v)$, $Q_P(v) = Q_P(-v)$ and $R_P(v) = R_P(-v)$.

Also, for any vertex $v$ of $P$ we have

(5.1) $$\mathrm{vol}_k P = \mathrm{vol}_k Q_P(v) + \mathrm{vol}_k N_P(v);$$

and

(5.2) $$\mathrm{vol}_k N_P(v) > 0.$$

By definition put $P_\lambda(v) = \mathrm{co}\{R(v), \pm \lambda v\}$. The idea behind this notation is that we will slightly scale one vertex of the projection of the cross-polytope and will write the first-order necessary condition for such an operation.

The next lemma is an immediate corollary of Lemma 1 from [RS58].

**Lemma 5.1.** *Let $v$ be a vertex of a centrally symmetric polytope $P$. Then $\mathrm{vol}_k P_\lambda(v)$ is a convex function of $\lambda$ and there exists $\varepsilon > 0$ such that*

(5.3) $$\mathrm{vol}_k P_\lambda(v) = \mathrm{vol}_k Q_P(v) + \lambda \, \mathrm{vol}_k N_P(v),$$

*for $\lambda \in [1, 1+\varepsilon]$.*

**Proof.**
By symmetry, we assume that $\lambda \geqslant 0$. All affine hyperplanes containing facets of $R_P(v)$ divide the line $\{tv \mid t \in \mathbb{R}\}$ into several intervals. As $\lambda$ increases, $\lambda v$ passed through these intervals. On each interval $\lambda v$ sees a certain set of facets of $P$ and hence, the function $\mathrm{vol}_k P_\lambda(v)$ is linear on each interval. As $\lambda v$ moves to the next interval, $\lambda v$ sees the larger sets of facets. Hence, the slope of the linear functions increase. This means that we have proven the convexity of the function.

For small enough positive $\varepsilon$ the interval $(v, (1+\varepsilon)v]$ belongs to one of the above mentioned intervals. It is easy to see that for each facet $F$ incident to the vertex $v$ the volume of $\mathrm{co}\{F, \lambda v, 0\}$ is $\lambda \mathrm{vol}_k \mathrm{co}\{F, 0\}$ for $\lambda \geqslant 1$. By the choice of $\varepsilon$ and by construction, we see that

$$\mathrm{vol}_k P_\lambda(v) = \mathrm{vol}_k Q_P(v) + \sum_{F \in \mathcal{F}(v), F \in \mathcal{F}(-v)} \mathrm{vol}_k(\mathrm{co}\{F, \lambda v, 0\}) = \mathrm{vol}_k Q_P(v) + \lambda \, \mathrm{vol}_k N_P(v).$$

This completes the proof.
$\square$

If the function $\mathrm{vol}_k P_\lambda(v)$ is linear in some open neighborhood of 1, then we can write a first-order necessary condition for maximizers and minimizers of (2.2). Moreover, this is the main observation in our proof of Theorem 1.5. Unfortunately, the function $\mathrm{vol}_k P_\lambda(v)$ is not linear in some open neighborhood of 1 in general, for example, see the corresponding polytopes for 3-dimensional cube in Figure 3. One can show that this function is linear in some neighborhood of 1 if and only if all facets incident to $v$ are simplices.

**Lemma 5.2.** *Let $\{\pm v_1, \ldots, \pm v_n\}$ be a set of pairwise distinct vectors and the vertex set of a centrally symmetric polytope $P$. Then*

$$\sum_{i \in [n]} \frac{\mathrm{vol}_k N(v_i)}{\mathrm{vol}_k P} \geqslant k.$$

*The bound is tight iff $P$ is a simplicial polytope.*

**Proof.**
Let $F$ be a facet of $P$ with $f$ vertices. Then the ratio $\frac{\mathrm{vol}_k(\mathrm{co}\{F,0\})}{\mathrm{vol}_k P}$ is a summand of exactly $f$ ratios $\frac{\mathrm{vol}_k N(v)}{\mathrm{vol}_k P}$. Since $f \geqslant k$, we have that

$$\sum_{i \in [n]} \frac{\mathrm{vol}_k N(v_i)}{\mathrm{vol}_k P} \geqslant k \sum_{F \in \mathcal{F}_P} \frac{\mathrm{vol}_k(\mathrm{co}\{F, 0\})}{\mathrm{vol}_k P} = k \frac{\mathrm{vol}_k P}{\mathrm{vol}_k P} = k,$$

where $\mathcal{F}_P$ is the set of all facets of $P$. Obviously, we have equality here only iff $P$ is a simplicial polytope.
$\square$

Now we are ready to prove a geometric necessary condition for (2.2).

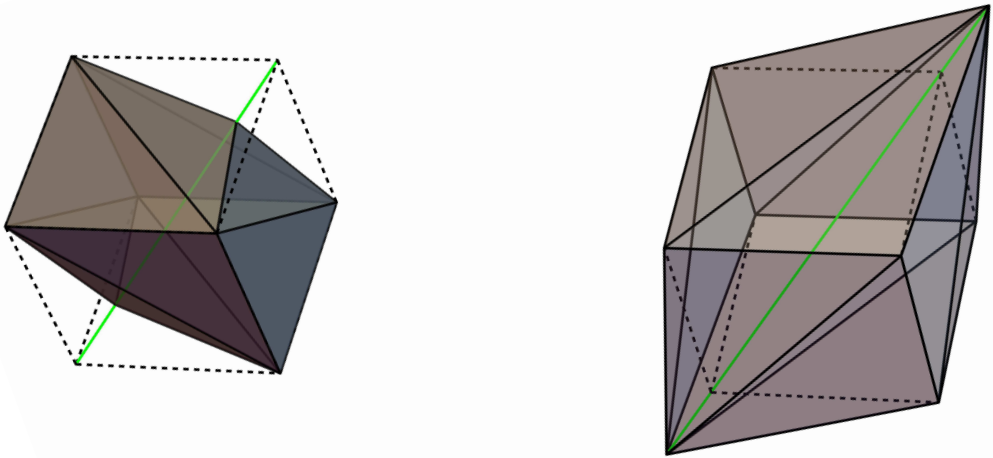

Figure 3. Scaling the corresponding opposite vertices of the cube along the green line, we obtain different combinatorial structure and linear coefficients.

**Theorem 5.3.** *Let a $S = \{v_1, \ldots, v_n\}$ be a tight frame such that a local maximum of (2.2) is attained on $S$. Then $\Diamond^n|S$ is a simplicial polytope and for every vertex $v$ of $\Diamond^n|S$ we have that $|v|^2 = \frac{\mathrm{vol}_k N_{(\Diamond^n|S)}(v)}{\mathrm{vol}_k(\Diamond^n|S)}$.*

**Proof.**
Denote $P = \Diamond^n|S$. Let $v_i \in S$ be a vertex of $P$. Let $\tilde{S}$ be the frame obtained from $S$ by substitution $v_i \to v_i + tv_i$. Substituting identities given by Lemma 3.3 and Lemma 5.1 into inequality (3.1), we get

$$1 + t|v_i|^2 + o(t) \geqslant \frac{\mathrm{vol}_k Q_P(v_i) + (1+t)\mathrm{vol}_k N_P(v_i)}{\mathrm{vol}_k P} = 1 + t\frac{\mathrm{vol}_k N_P(v_i)}{\mathrm{vol}_k P},$$

for a sufficiently small positive $t$. From this, we obtain that $|v_i|^2 \geqslant \frac{\mathrm{vol}_k N_P(v_i)}{\mathrm{vol}_k P}$. By Lemma 4.4, we have that all non-zero elements of $\{\pm v_1, \ldots, \pm v_n\}$ are pairwise distinct vertices of $P$. Therefore, by (4.2), we have

$$k = \sum_{i \in [n] | v_i \neq 0} |v_i|^2 \geqslant \sum_{i \in [n] | v_i \neq 0} \frac{\mathrm{vol}_k N_P(v_i)}{\mathrm{vol}_k P}.$$

And finally, by Lemma 5.2, the right-hand side is at least $k$ and it is $k$ if and only if $P$ is a simplicial polytope. Therefore, all inequalities are tight. This completes the proof. $\square$

Theorem 5.3 implies that if $S \in \Omega(n,k)$ is a local maximum of (2.2) then for an arbitrary sufficiently small perturbation of $S$, the combinatorial structure of the associated polytope remains the same (for a sufficiently small perturbation, the images of zero vectors of $S$ remain in the interior and do not affect the volume). Hence, the volume of an associated polytope is a differentiable function in a suffiently small neighborhood of a local maximum. We can use all operations described in Section 3, which, along with Lemma 3.2, give us some geometric restriction for a tight frame to be a local maximum.

# 6. Proofs of Theorem 1.5 and Theorem 1.2

**Proof of Theorem 1.5.** The main observation in this proof is that polygons are simplicial polytopes. Therefore, we can use identity (5.3) for the projection. Let $S = \{v_1, \ldots, v_n\}$ be a tight frame such that $\text{vol}_2\,(\Diamond^n | S)$ is the global minimum of (2.2). Denote $P = \Diamond^n | S$ and let $v$ be a vertex of $P$. Clearly, $v \in \{\pm v_1, \ldots, \pm v_n\}$, and, by Corollary 4.5, we know that $\{\pm v_1, \ldots, \pm v_n\}$ are vertices of $P$. Let exactly $d$ pairs $\pm v_1, \ldots, \pm v_d$ coincide with $\pm v$. Let $\tilde{S}$ be the frame obtained from $S$ by substitution $v_i \to (1+t)v_i$, for all $i \in [d]$ and $t \in (-\varepsilon, \varepsilon)$ for a sufficiently small $\varepsilon$ such that we do not change a combinatorial type of $P$.

By identity (4.1) and identity
$$A_{\tilde{S}} = A_S - dv \otimes v + d(1+t)^2 v \otimes v = I_2 + (2t + t^2) dv \otimes v,$$
we obtain that $\det A_{\tilde{S}} = 1 + (2t + t^2)d|v|^2$. By identity (5.3), the ratio in the left-hand side of inequality (3.1) equals
$$1 + t\frac{\text{vol}_2\,N_P(v)}{\text{vol}_2\,P(v)},$$
for positive $t < \varepsilon$. By the choice of $\varepsilon$, the same identity holds on the whole interval $(-\varepsilon, \varepsilon)$. By Lemma 3.2, we obtain
$$1 + t\frac{\text{vol}_2\,N_P(v)}{\text{vol}_2\,P(v)} \geqslant \sqrt{1 + (2t + t^2)d|v|^2} \quad \text{for all } t \in (-\varepsilon, \varepsilon).$$
Consider the Taylor expansion of the right-hand side. Since the inequality holds for all $t \in (-\varepsilon, \varepsilon)$, the linear part in $t$ on the both sides of the inequality coincide. Since there is no quadratic part in the left-hand side, the coefficient of $t^2$ in the right-hand side must be non-positive. This coefficient is
$$C = \frac{(d|v|^2)^2 - d|v|^2}{2}.$$
Since $S$ is a frame and by Lemma 2.3, we have that $d|v|^2 \in (0, 1]$. Therefore, $C \leqslant 0$ if and only if $d|v|^2 = 1$. Again, by Lemma 2.3, this implies that $\langle v_j, v \rangle = 0$ for $j > d$, and, by (4.2), $P$ has only two pairs of centrally symmetric vertices. Or, equivalently, $P$ is a rhombus. The optimal constant in this case was found in Lemma 4.7.
□

Actually, we used only that $S$ is a local minimum to prove that $\Diamond^n | S$ is a rhombus in Theorem 1.5. This means that all local minima are attained on rhombi. But the maximization problem is more sophisticated.

**Proof of Theorem 1.2.** The main idea is to use substitution $v \to 0$ for an appropriately chosen vector $v$ of a maximizer. Such an operation is not continuous, and, as will be shown in the next section, we cannot transform an arbitrary tight frame to obtain the global maximum while increasing the volume of the generated projection of the cross-polytope in a monotonic way.

So, let $S = \{v_1, \ldots, v_n\}$ be a tight frame such that the global maximum of (2.2) is attained on it, and denote $P = \Diamond^n | S$. By Lemma 4.4, we assume that $\{\pm v_1, \ldots, \pm v_n\}$ are pairwise distinct vertices of $P$. Choose a vertex $v$ of $P$ such that $|v| < 1$ (if $|v_i| = 1$ for all $i \in [n]$ then, by Lemma 2.3, $n = k$ and the vectors of $S$ form an orthonormal basis, and $P$ is the standard cross-polytope). Let $t \in (0, 1)$ be such that $tv$ is on the boundary of the set $R_P(v)$ (if $t = 0$ then all other vectors in $S$ are orthogonal to $v$, and hence, $|v| = 1$). Let us prove that $t$ satisfies the following inequality

(6.1) $$\frac{1}{2} > \frac{x}{1+x} \geqslant t \geqslant \frac{x}{\sqrt{n-1}\sqrt{1-x^2}},$$

where $x = \sqrt{1 - |v|^2} \in (0, 1)$.

First, we estimate $t$ from above. We assume $v = v_1$. Let $\tilde{S}$ be the frame obtained from $S$ by substitution $v_1 \to 0$. Then, by (4.1) and by Lemma 3.2, we obtain

$$\sqrt{1 - |v^2|} \geqslant \frac{\operatorname{vol}_k\left(\Diamond^n|\tilde{S}\right)}{\operatorname{vol}_k P}. \tag{6.2}$$

But $\Diamond^n|\tilde{S}$ is precisely $R_P(v)$. The rest of $v$ in $P$ is the union of the two internally disjoint sets
  (1) the belt $Q_P(v)$, which, by Theorem 5.3, has the volume $(1 - |v|^2)\operatorname{vol}_k P$;
  (2) the intersection of the rest and the star of the vertex $v$; this set contains $tv$, and, therefore, its volume is at least $t\operatorname{vol}_k N_P(v)$, which is, by Theorem 5.3, $t|v|^2 \operatorname{vol}_k P$.

Returning to inequality (6.2), we obtain

$$\sqrt{1 - |v|^2} \geqslant (1 - |v|^2) + t|v|^2 \Leftrightarrow \frac{\sqrt{1 - |v|^2} - (1 - |v|^2)}{|v|^2} \geqslant t \Leftrightarrow \frac{x}{1 + x} \geqslant t. \tag{6.3}$$

Second, we estimate $t$ from below. Let $H$ be the supporting hyperplane to $R_P(v)$ at $tv$ and $\ell$ be its orthogonal complement. Then, by Lemma 3.4, the projections $\{v'_1, \ldots, v'_n\}$ of vectors of $S$ onto $\ell$ form a tight frame in $\ell$, hence $\sum_{i \in [n]} |v'_i|^2 = 1$. On the other hand, the projection of $v$ has the length at most $|v|$ and the projections of all others have the length at most $t|v|$. This means that we have proven inequality $|v|^2 + (n-1)t^2|v|^2 \geqslant 1$, or, equivalently,

$$t \geqslant \frac{x}{\sqrt{n-1}\sqrt{1-x^2}}.$$

This completes the proof of inequality (6.1).

Let us complete the proof of the theorem using inequality (6.1). As a straightforward consequence of (6.1), we obtain the following chain of equivalent inequalities:

$$\frac{x}{1+x} \geqslant \frac{x}{\sqrt{n-1}\sqrt{1-x^2}} \Leftrightarrow \frac{1-x}{1+x} \geqslant \frac{1}{n-1} \Leftrightarrow 1 - \frac{2}{n} \geqslant x.$$

By identity (4.2), $v \in S$ can be chosen such that $x \geqslant \sqrt{1 - \frac{k}{n}}$. Substituting this in the last inequality, we obtain

$$1 - \frac{4}{n} + \frac{4}{n^2} \geqslant 1 - \frac{k}{n} \Leftrightarrow k \geqslant 4 - \frac{4}{n}.$$

The latter does not hold neither when $k = 2$ and $n \geqslant 3$ nor when $k = 3$ and $n > 4$. For $k = 2$, this means that the global maximum of (2.2) is attained on a tight frame with exactly $n = 2$ non-zero vectors. By Lemma 2.3, this yields that these two vectors form an orthonormal basis of their linear hall. And for $k = 3$, we have that the global maximum of (2.2) is attained on a tight frame with either $n = 3$ or $n = 4$ non-zero vectors. As mentioned in the Introduction, the hyperplane case (in particular, $k = 3$ and $n = 4$) was proven by K. Ball in [Bal95]. Therefore, a global maximizer for $k = 3$ has exactly 3 non-zero vectors, and again, by Lemma 2.3, these three vectors form an orthonormal basis of their linear hall.
$\square$

## 7. Local maxima

**Lemma 7.1.** *Let $W = \{w_1, w_2, w_3\}$ be a tight frame in $\mathbb{R}^2$ such that $\Diamond^3|W$ is a regular hexagon. Then $W$ is a local maximum of (2.2).*

**Proof.**
Let $\Omega \subset \Omega(3, 2)$ be the set of all tight frames $S = \{v_1, v_2, v_3\}$ in $\mathbb{R}^2$ such that all three vectors $v_1, v_2, v_3$ are on the boundary of $\Diamond^3|S$. It is a closed set. Hence, there exists a global maximum of (2.2) on $\Omega$. We will show that $W$ is the unique global maximizer on $\Omega$ upto orthogonal transformations. Since some neighborhood of $W$ in $\Omega(3, 2)$ belongs to the interior of $\Omega$, this proves that $W$ is a local maximum of (2.2) on $\Omega(3, 2)$.

Let the global maximum on $\Omega$ be attained on a tight frame $S = \{v_1, v_2, v_3\}$.

We begin with excluding the degenerate case when the associated polytope $\Diamond^3 | S$ is a parallelogram, we assume $v_2 \in [v_1, v_3]$. Let us show that the vector $v_2$ is perpendicular to the segment $[v_1, v_3]$. Assume the contrary, then the circle with center at the origin and the radius $|v_2|$ intersects $[v_1, v_3]$ at point $v_2$ and the segment $[v_1, v_3]$ is not tangent to the circle. Therefore, slightly rotating $v_2$ around the origin, we get a point $v_2'$ outside $\Diamond^3 | S$. Consider a new frame $\tilde{S} = \{v_1, v_2', v_3\}$. Clearly, $\text{vol}_2 \left( \Diamond^3 | \tilde{S} \right) > \text{vol}_2 \left( \Diamond^3 | S \right)$. However, since $\text{tr} A_{\tilde{S}} = \text{tr} A_S = 2$ and by the inequality of arithmetic and geometric means, we have that $\det A_{\tilde{S}} < \det A_S = 1$. Therefore, by Lemma 3.2, $S$ is not a local maximum. We proved that $v_2 \perp [v_1, v_3]$. Choosing $v_2$ as the direction of the first vector of an orthonormal basis of the plane and using assertion (4) of Lemma 2.3, we obtain that

$$S = \begin{pmatrix} 1/\sqrt{3}, & 1/\sqrt{3}, & 1/\sqrt{3} \\ \pm 1/\sqrt{2}, & 0, & \mp 1/\sqrt{2} \end{pmatrix}.$$

Clearly, $\text{vol}_k \left( \Diamond^3 | S \right) < \text{vol}_k \left( \Diamond^3 | W \right)$.

Now, let $\Diamond^3 | S$ be a hexagon, see notation in Figure 4.

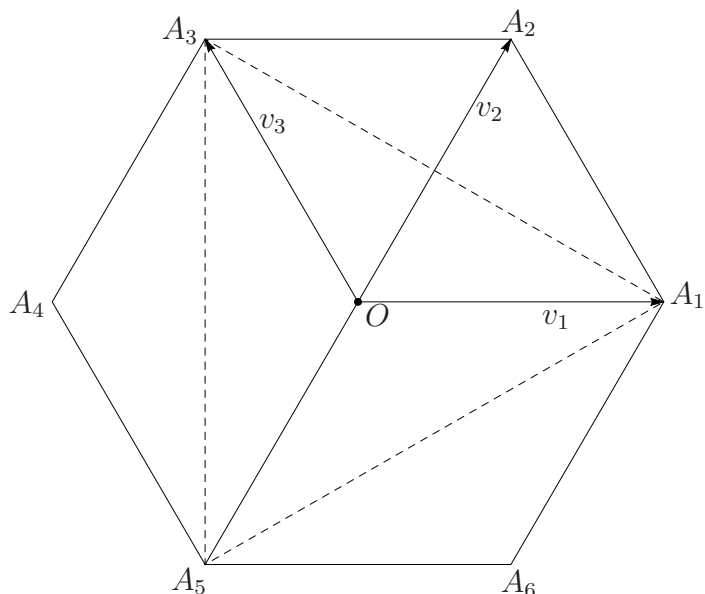

FIGURE 4. Illustration for Lemma 7.1.

In this case $S$ is in the interior of $\Omega$. Lemma 3.2 yields two observations.

(1) We claim that $v_2 \perp [v_1, v_3]$. We slightly rotate $v_2$ around the origin. We only change triangles $A_3 A_2 A_1$ and $A_6 A_5 A_4$. If $v_2$ is not perpendicular to $[v_1, v_3]$, then rotating in the right direction we increase the altitude through vertex $A_2$ and perpendicular to $A_1 A_3$ in $A_3 A_2 A_1$ and the altitude through vertex $A_5$ and perpendicular to $A_4 A_6$ in $A_4 A_5 A_6$, hence, we increase the area of the polygon $\Diamond^3 | S$. Again, by the inequality of arithmetic and geometric means, we decrease the determinant rotating $v_2$. This contradicts Lemma 3.2.

(2) We claim that the triangles $OA_1 A_2$ and $OA_2 A_3$ have the equal area. We use the substitution $v_1 \to v_1 + t v_3, v_3 \to v_3 - t v_1$. By Lemma 3.3, the linear part in the Taylor expansion of the determinant in the right-hand side of (3.1) vanishes. The area of the polygone $\Diamond^3 | S$ is

$$\det(v_1, v_2) + \det(v_2, v_3) + \det(v_3, -v_1).$$

Hence, the change of the area up to the terms of the first order equals

$$t(\det(v_3, v_2) - \det(v_2, v_1)) = -2t \left( \text{Area} \triangle OA_2 A_3 - \text{Area} \triangle OA_1 A_2 \right).$$

By Lemma 3.2, it has to vanish as well.

Since the areas of triangles $OA_1A_2$ and $OA_2A_3$ are equal, $OA_2$ is a median of the triangle $OA_1A_3$. Therefore, $A_5A_2$ is a median and an altitude of the triangle $A_1A_3A_5$, hence it is an isosceles triangle. By the symmetry, the triangles $A_1A_3A_5$ and $A_2A_4A_6$ are equilateral and $O$ is their common center. We conclude that $\Diamond^3|S$ is a regular hexagon. This completes the proof. □

Again, by Lemma 4.4, we can add zero vectors to $W$ and a new tight frame remains a local maximum of (2.2). Moreover, for a given local maximum for $\Omega(n,k)$ we can construct a new local maximum for $\Omega(n+1, k+1)$.

**Lemma 7.2.** *Let $S \subset \Omega(n,k)$ give a local maximum of (2.2). Consider $S_1$ obtained from $S$ as follows*

- *considering the standard embedding $\mathbb{R}^k \subset \mathbb{R}^{k+1}$, we extend all vectors of $S$ by adding zero in the last $(k+1)$-st coordinate;*
- *add $e_{k+1}$.*

*Then $S_1 \subset \Omega(n+1, k+1)$ and it gives a local maximum of (2.2) for $\Omega(n+1, k+1)$.*

**Proof.**
By Theorem 5.3, $\Diamond^n|S$ is a simplicial polytope. Hence, we can subdivide it into non-overlapping simplices of the type $\mathrm{co}\{F, 0\}$, where $F$ is a facet of $\Diamond^n|S$. Then $S_1$ is a simplicial polytope as well, and simplices $\mathrm{co}\{F, 0, \pm e_{k+1}\}$ subdivide $\Diamond^{n+1}|S_1$. Let $S'_1 = \{v_1, \ldots, v_{n+1}\} \subset \Omega(n+1, k+1)$ be a sufficiently small perturbation of $S_1$ (where $v_{n+1}$ is a perturbation of $e_{k+1}$) and $H^\perp$ be the orthogonal complement to $v_{n+1}$ in $\mathbb{R}^{k+1}$. Denote by $P_\perp$ the orthogonal projector onto $H^\perp$. Then, by linearity and since $\Diamond^{n+1}|S'_1$ and $\Diamond^{n+1}|S_1$ have the same combinatorial structure, we have

$$\mathrm{vol}_{k+1}\left(\Diamond^{n+1}|S'_1\right) = \sum \left(\mathrm{vol}_{k+1}\,\mathrm{co}\{F, 0, v_{n+1}\} + \mathrm{vol}_{k+1}\,\mathrm{co}\{F, 0, -v_{n+1}\}\right) =$$

$$\sum \left(\mathrm{vol}_{k+1}\,\mathrm{co}\{P_\perp F, 0, v_{n+1}\} + \mathrm{vol}_{k+1}\,\mathrm{co}\{P_\perp F, 0, -v_{n+1}\}\right) = \mathrm{vol}_{k+1}\left(\Diamond^{n+1}|S_2\right),$$

where the summations are over all $(k-1)$-dimensional faces of $\Diamond^{n+1}|S'_1$ generated by the vectors $\{\pm v_1, \ldots, \pm v_n\}$ and $S_2 = \{P_\perp v_1, \ldots, P_\perp v_n, v_{n+1}\}$.

By Lemma 3.4, $S_3 = \{P_\perp v_1, \ldots, P_\perp v_n\}$ form a tight frame in $H^\perp$ and it can be considered as a sufficiently small perturbation of $S$. Therefore, by the choice of $S$, we have

$$\mathrm{vol}_{k+1}\left(\Diamond^{n+1}|S'_1\right) = \mathrm{vol}_{k+1}\left(\Diamond^{n+1}|S_2\right) = \frac{2}{k+1}|v_{n+1}|\,\mathrm{vol}_k\left(\Diamond^n|S_3\right) \leqslant \frac{2}{k+1}\,\mathrm{vol}_k\left(\Diamond^n|S\right) = \mathrm{vol}_k\left(\Diamond^{n+1}|S_1\right).$$

That is, $S_1$ gives a local maximum. □

Summarizing, we have just proved.

**Theorem 7.3.** *For any $n > k \geqslant 2$ there exists a local and not the global maximum of (2.2).*

Therefore, it is not sufficient to use only a necessary condition of local maximum (for example, Theorem 5.3) in order to prove Conjecture 1.1. Our approach in a lower dimensional case is to map one vector to the origin, and it works for $k = 2, 3$. We do not know whether our approach works in a higher dimensional case. The following questions remain open: can we delete one (or even several) vertex, transform to a new tight frame (using the scheme $S \to \tilde{S} \to S'$) and increase the volume for an arbitrary polytope in $\mathbb{R}^n$? Is it true that the maximum volume of the projection of the cross-polytope with $2\ell$ vertices is a decreasing function of $\ell$?

## 8. Appendix

**Sketch of the proof of Lemma 3.3.**
Recall that the cross product of $k-1$ vectors $\{x_1, \ldots, x_{k-1}\}$ of $\mathbb{R}^k$ is the vector $x$ defined by

$$\langle x, y \rangle = \det(x_1, \ldots, x_{k-1}, y) \quad \text{for all} \quad y \in \mathbb{R}^k.$$

For an ordered $(k-1)$-tuple $L = \{i_1, \ldots, i_{k-1}\} \in \binom{[n]}{k-1}$ and a frame $S = \{v_1, \ldots, v_n\}$, we use $[v_L]$ to denote the cross product of $v_{i_1}, \ldots, v_{i_{k-1}}$.

We claim the following property of the tight frames.

Let $S = \{v_1, \ldots, v_n\}$ be a tight frame in $\mathbb{R}^k$. Then vectors $\{[v_L]\}_{L \in \binom{[n]}{k-1}}$ form a tight frame in $\mathbb{R}^k$.

We use $\Lambda^k(\mathbb{R}^n)$ to denote the space of $k$-forms on $\mathbb{R}^n$. By assertion (2) of Lemma 2.3, there exist an orthonormal basis $\{f_i\}_1^n$ of $\mathbb{R}^n$ such that $v_i$ is the orthogonal projection of $f_i$ onto $\mathbb{R}^k$, for any $i \in [n]$. Then the $(k-1)$-form $v_{i_1} \wedge \cdots \wedge v_{i_{k-1}}$ is the orthogonal projections of the $(k-1)$-form $f_{i_1} \wedge \cdots \wedge f_{i_{k-1}}$ onto $\Lambda^{k-1}(\mathbb{R}^k) \subset \Lambda^{k-1}(\mathbb{R}^n)$, for any ordered $(k-1)$-tuple $L = \{i_1, \ldots, i_{k-1}\} \in \binom{[n]}{k-1}$. By Lemma 2.3 and since the $(k-1)$-forms $\{f_{i_1} \wedge \cdots \wedge f_{i_{k-1}}\}_{\{i_1, \ldots, i_{k-1}\} \in \binom{[n]}{k-1}}$ form an orthonormal basis of $\Lambda^{k-1}(\mathbb{R}^n)$, we have that $(k-1)$-forms $\{v_{i_1} \wedge \cdots \wedge v_{i_{k-1}}\}_{\{i_1, \ldots, i_{k-1}\} \in \binom{[n]}{k-1}}$ form a tight frame in $\Lambda^{k-1}(\mathbb{R}^k)$. Finally, the Hodge star operator maps $v_{i_1} \wedge \cdots \wedge v_{i_{k-1}}$ to the cross product of vectors $v_{i_1}, \ldots, v_{i_{k-1}}$. Since the Hodge star is an isometry, the cross products $\{[v_L]\}_{L \in \binom{[n]}{k-1}}$ form a tight frame. The claim is proven.

By linearity of the determinant, it is enough to prove the lemma for $S' = \{v_1 + tx, v_2, \ldots, v_n\}$. Denote $v_1' = v_1 + tx$ and $v_i' = v_i$, for $2 \leq i \leq n$.

By the Cauchy-Binet formula, we have

$$(8.1) \quad \det A_{S'} = \det\left(\sum_1^n v_i' \otimes v_i'\right) = \sum_{Q \in \binom{[n]}{k}} \det\left(\sum_{i \in Q} v_i' \otimes v_i'\right).$$

By the properties of the Gram matrix, we have

$$\det\left(\sum_1^k v_{i_1}' \otimes v_{i_k}'\right) = \left(\det(v_{i_1}', \ldots, v_{i_k}')\right)^2.$$

By this, by the definition of cross product and by (8.1), we obtain

$$\det A_{S'} = 1 + 2t \sum_{Q \in \binom{[n]}{k}, 1 \in Q} \langle v_1, [v_{Q\setminus 1}]\rangle \langle [v_{Q\setminus 1}], x\rangle + o(t).$$

Since $\langle v_1, [v_J]\rangle = 0$ for any $J \in \binom{[n]}{k-1}$ such that $1 \in J$, we have that the linear term of the Taylor expansion of $\det A_{S'}$ equals

$$2t \sum_{Q \in \binom{[n]}{k}, 1 \in Q} \langle v_1, [v_{Q\setminus 1}]\rangle \langle [v_{Q\setminus 1}], x\rangle = 2t \sum_{L \in \binom{[n]}{k-1}} \langle v_i, [v_L]\rangle \langle [v_L], x\rangle.$$

Since $\{[v_L]\}_{L \in \binom{[n]}{k-1}}$ form a tight frame in $\mathbb{R}^k$, we have that

$$\sum_{L \in \binom{[n]}{k-1}} \langle v_i, [v_L]\rangle \langle [v_L], x\rangle = \langle v_i, x\rangle.$$

Therefore,

$$\sqrt{\det A_{S'}} = \sqrt{\det A_S + 2t\langle v_i, x\rangle + o(t)} = 1 + t\langle v_i, x\rangle + o(t).$$

□

## References

[AHK19] Arseniy Akopyan, Alfredo Hubard, and Roman Karasev. Lower and upper bounds for the waists of different spaces. *Topological Methods in Nonlinear Analysis*, 53(2):457–490, 2019.

[Bal89] Keith Ball. Volumes of sections of cubes and related problems. *Geometric aspects of functional analysis*, pages 251–260, 1989.

[Bal95] Keith Ball. Mahler's conjecture and wavelets. *Discrete & Computational Geometry*, 13(1):271–277, 1995.

[Bar98]    Franck Barthe. On a reverse form of the Brascamp-Lieb inequality. *Inventiones mathematicae*, 134(2):335–361, 1998.

[BN02]     Franck Barthe and A Naor. Hyperplane projections of the unit ball of $\ell_p^n$. *Discrete & Computational Geometry*, 27(2):215–226, 2002.

[Fil88]    P. Filliman. The largest projections of regular polytopes. *Israel Journal of Mathematics*, 64(2):207–228, Jun 1988.

[Iva17]    Grigory Ivanov. On the Volume of the John–Löwner Ellipsoid. *Discrete & Computational Geometry*, pages 1–5, 2017.

[Iva18]    Grigory Ivanov. Tight frames and related geometric problems. *arXiv e-prints*, Apr 2018.

[LYZ$^+$04] Erwin Lutwak, Deane Yang, Gaoyong Zhang, et al. Volume inequalities for subspaces of $L_p$. *Journal of Differential Geometry*, 68(1):159–184, 2004.

[MP88]     Mathieu Meyer and Alain Pajor. Sections of the unit ball of $\ell_p^n$. *Journal of functional analysis*, 80(1):109–123, 1988.

[RS58]     CA Rogers and GC Shephard. Some extremal problems for convex bodies. *Mathematika*, 5(2):93–102, 1958.

[Vaa79]    Jeffrey Vaaler. A geometric inequality with applications to linear forms. *Pacific Journal of Mathematics*, 83(2):543–553, 1979.

[Zon06]    Chuanming Zong. *The Cube – A Window to Convex and Discrete Geometry*, volume 168. Cambridge University Press, 2006.